\documentclass[12pt]{article}
\usepackage{euscript,amsmath, amssymb, amsfonts}
\usepackage{color}

\newcommand{\HH}{\mathcal H }

\newcommand{\G}{{\mathbb C}[I_2(2m+1)]}
\newcommand{\spur}{sp}

\newcommand{\HG}{$\mathcal H$}

\newcommand{\m}{\mathfrak {s} }
\newcommand{\Lo}{\mu L_0 }

\newcommand{\be}{\begin{equation}}
\newcommand{\ee}{\end{equation}}
\newcommand{\bee}{\begin{eqnarray}}
\newcommand{\eee}{\end{eqnarray}}
\newcommand\nn{\nonumber \\}

\newcommand\defeq{:=}

\newcounter{theorem}
\newcommand{\theorem}{\par\refstepcounter{theorem}
           {\bf Theorem \arabic{section}.%
           \arabic{theorem}. }}
\renewcommand\thetheorem{\thesection.%
   \arabic{theorem}}
\makeatletter \@addtoreset{theorem}{section}

\newcounter{corollary}
\newcommand{\corollary}{\par\refstepcounter{corollary}
           {\bf Corollary \arabic{section}.%
           \arabic{theorem}. }}

\makeatletter \@addtoreset{corollary}{section}

\newcounter{lemma}

\makeatletter \@addtoreset{lemma}{section}

\newcounter{proposition}
\newcommand{\proposition}{\par\refstepcounter{theorem}
           {\bf Proposition \arabic{section}.%
           \arabic{theorem}. }}
\renewcommand\theproposition{\thesection.%
  \arabic{theorem}}
\makeatletter \@addtoreset{proposition}{section}

\newcounter{conjecture}
\newcommand{\conjecture}{\par\refstepcounter{theorem}
           {\bf Conjecture \arabic{section}.%
           \arabic{theorem}. }}

\makeatletter \@addtoreset{conjecture}{section}

\newcounter{remark}
\newcommand{\remark}{\par\refstepcounter{theorem}
           {\bf Remark \arabic{section}.%
           \arabic{theorem}. }}
\renewcommand\theremark{\thesection.%
     \arabic{theorem}}
\makeatletter \@addtoreset{remark}{section}

\newcounter{definition}
\newcommand{\definition}{\par\refstepcounter{theorem}
           {\bf Definition \arabic{section}.%
           \arabic{theorem}. }}
\renewcommand\thedefinition{\thesection.%
    \arabic{theorem}}
\makeatletter \@addtoreset{definition}{section}

\newenvironment{proof}[1][Proof]{\noindent\textsf{#1.\ }}
{\hfill {\small $\square$}}

\makeatletter \@addtoreset{equation}{section}
\def\theequation{\thesection.\arabic{equation}}

\begin{document}


\sloppy \title
 {
Ideals generated by traces or by supertraces in the symplectic reflection
algebra $H_{1,\nu}(I_2(2m+1))$
}

\author
 {
 S.E. Konstein%
\thanks{ I.E. Tamm Department of Theoretical Physics,
          P.N. Lebedev Physical Institute, RAS
          119991, Leninsky prosp., 53, Moscow, Russia}
\thanks{E-mail: konstein@lpi.ru}
 ,
 I.V. Tyutin$^*$%
\thanks{Tomsk State Pedagogical University, Tomsk, Russia}
\thanks{E-mail: tyutin@lpi.ru}     }

\date{
}

\maketitle
\thispagestyle{empty}

\begin{abstract}

For each complex number $\nu$,
an associative symplectic reflection algebra $\mathcal H:= H_{1,\nu}(I_2(2m+1))$,
based on the group generated by root system $I_2(2m+1)$,
has an $m$-dimensional space of traces and an $(m+1)$-dimensional space of supertraces.
A (super)trace $\spur$ is said to be degenerate if the corresponding bilinear
(super)symmetric form $B_{sp}(x,y)=sp(xy)$ is degenerate. We find all values of the parameter
 $\nu$ for which either the space of traces contains a degenerate nonzero trace or
the space of supertraces contains a degenerate nonzero supertrace and, as a consequence,
the algebra $\mathcal H$ has a two-sided ideal of null-vectors.

The analogous results for the case $H_{1,\nu_1, \nu_2}(I_2(2m))$ are also presented.
\end{abstract}


\section{Introduction}

\subsection{Definitions}

Let $\mathcal A$ be an associative $\mathbb Z_2$-graded algebra with unit
and with a parity $\varepsilon$. All expressions of linear algebra are given
for homogenous elements only and are supposed to be extended to inhomogeneous elements via
linearity.

A linear complex-valued function $tr$ on $\mathcal A$ is called a \emph{trace}
if $tr(fg-gf)=0$ for all $f,g\in \mathcal A$.
A linear complex-valued function $str$ on $\mathcal A$ is called a \emph{supertrace}
if ${str(fg-(-1)^{\varepsilon(f)\varepsilon(g)}gf)=0}$ for all $f,g\in \mathcal A$.
These two definitions can be unified as follows.

Let ${\varkappa=\pm 1}$.
A linear complex-valued function  $\spur$ on $\mathcal A$ is called \emph{$\varkappa$-trace}
if ${\spur(fg-\varkappa^{\varepsilon(f)\varepsilon(g)}gf)=0}$ for all $f,g\in \mathcal A$.

The element $K\in \mathcal A$ is called \emph {Klein operator} if $K^2=1$, $\varepsilon(K)=0$,
and $Kf=(-1)^{\varepsilon(f)}fK$ for any $f\in \mathcal A$. If the algebra $\mathcal A$ contains
Klein operator $K$, then the linear function $f\mapsto tr(fK^{\varepsilon(f)+1})$ is a supertrace,
and the linear function $f\mapsto str(fK^{\varepsilon(f)+1})$ is a trace.

Each nonzero (super)trace $\spur$ defines the nonzero (super)symmetric bilinear form
 $B_{\spur}(f,g):=\spur(fg)$.

If this bilinear form is degenerate, then the set of its null-vectors is a proper
ideal in $\mathcal A$. We say that the (super)trace $\spur$ is \emph{degenerate}
if the bilinear form $B_{\spur}$ is degenerate.

\subsection{The goal and structure of the paper}

In a number of papers the simplicity (or, alternatively, existence of ideals) of
Symplectic Reflection Algebras or. briefly, SRA (for definition, see \cite{sra})
was investigated, see, e.g.,  \cite{BG}, \cite{IL}. In particular, it is shown that
all SRA with zero parameters of deformation
are simple (see \cite{pass}, \cite{BG}).

It follows from \cite{KT} and \cite{stek} that an associative algebra of observables of
the Calogero model with harmonic term in the potential and with coupling constant $\nu$
based on the root system $I_2(2m+1)$ (this algebra is SRA $H_{1,\nu}(I_2(2m+1))$)
has an $m$-dimensional space of the traces and
an $(m+1)$-dimensionsl space of supertraces.

We say that the parameter $\nu$ is \emph{singular}, if the algebra $H_{1,\nu}(I_2(n))$ has
a degenerate trace or supertrace.

The goal of this paper is to find all singular values of $\nu$
for the algebras $\mathcal H := H_{1,\nu}(I_2(n))$ and to find
corresponding degenerate traces and supertraces  for $n$ odd, ${n=2m+1}$;
the result is formulated in Theorem \ref{th2}.

The consideration
generalizes \cite{K2}, where ${H_{1,\nu}(A_2) \simeq   H_{1,\nu}(I_2(3))}$ was considered.

The case of $n$ even is considered in \cite{tmf} and reproduced here in Appendix.
It is shown there that the set of singular values of $\nu$-s consists of
four families of parallel complex lines on the complex plane $(\nu_1,\,\nu_2)$.
The results for $H_{1,\nu}(I_2(2m))$ generalize
the simplest case
$H_{1,\nu_1, \nu_2}(I_2(2)) \simeq H_{1,\nu_1}(A_1)\otimes H_{1,\nu_2}(A_1)$;
the singular values of $\nu$ and ideals in $H_{1,\nu}(A_1)$
were found in \cite{V}.

In Sections \ref{i2}-\ref{genfuni2} we recall the necessary definitions
and prove the preliminary facts.
In Section \ref{degcond}, we show
that if the $\varkappa$-trace is degenerate, then its generating functions
are integer.
In Section \ref{genfun}, we derive the equations for the generating functions
of the $\varkappa$-trace and solve these equations. The solutions are meromorphic
functions of their parameter $t$ for every value of $\nu$, except degenerate values,
which are found in Section \ref{spec}.

\section{The group $I_2(n)$}\label{i2}

\definition
The group $I_2(n)$ is a finite subgroup of the orthogonal group $O(2,\mathbb R)$ generated
by the root system $I_2(n)$.

In this subsection we consider ${\mathbb C}$ instead of ${\mathbb R}^2$ for convenience.

The root system $I_2(n)$ consists of $2n$ vectors $v_k=\exp(\frac {i \pi k} n)$, where
$k=0,1,...,2n-1$. The group $I_2(n)$
has $2n$ elements: $n$ reflections $R_k$, acting on $z \in{\mathbb C}$
as follows:
\bee\label{Rk}
&& R_k z   = -z^* v_k^2 R_k, \ \text{ for } k=0,\,1,\,\dots,\,n-1,
\eee
and $n$ rotations, which are the elements of the form $S_k:=R_k R_0$. The element $S_0$
is the unit in the group $I_2(n)$.
It is easy to see from the formula (\ref{Rk}) that these elements satisfy the relations
\be\label{RR}     \nonumber
R_k R_l = S_{k-l},\qquad
S_k S_l = S_{k+l},\qquad
R_k S_l = R_{k-l},\qquad
S_k R_l = R_{k+l}.
\ee
Obviously, if $n$ is even, then the reflections $R_{2k}$, where $k=0,2,...,\frac n 2 - 1$,
constitute one
conjugacy class, and the $R_{2k+1}$ constitute another class. If $n$ is odd,
then all the reflections $R_k$
are in the same conjugacy class.

The rotations $S_k$ and $S_l$ constitute a conjugacy class if
$k+l=n$.

Let
\be\label{defg}
G:={\mathbb C}[I_2(n)]
\ee
be the group algebra of the group $I_2(n)$.
In $G$, it is convenient to introduce the following basis
\bee\label{LQbas}
&& L_p:=\frac 1 n \sum_{k=0}^{n-1} \lambda^{kp} R_k,\qquad
Q_p:=\frac 1 n \sum_{k=0}^{n-1} \lambda^{-kp} S_k,\\
\label{lambda}\nonumber
&& \mbox{ where }\lambda=\exp\left(\frac{2\pi i} n\right).
\eee

In what follows, we consider only the case of $n$ odd, $n=2m+1$.

The result for $n$ even from \cite{tmf} is reproduced in Appendix.
The main differences between odd and even $n$ are as follows:

an even) the algebra $H_{1,\nu}(I_2(2m))$ depends on two complex parameters $\nu$;
the algebra $H_{1,\nu}(I_2(2m))$ contains the Klein operator and so
the space of traces and the space of supertraces are isomorphic;

an odd) the algebra $H_{1,\nu}(I_2(2m+1))$ depends on one complex parameter $\nu$,
the space of traces is an $m$-dimensional, the space of supertraces is
an $(m+1)$-dimensional.


\section{Symplectic reflection algebra
$H_{1,\nu}(I_2(2m+1))$}\label{srai2}

Let $n$ be odd, $n=2m+1$, and $\lambda=\exp \left(\frac {2\pi i} n\right)$.

\definition
\emph{The symplectic reflection algebra
$\mathcal H :=H_{1,\nu}(I_2(2m+1))$} is the associative algebra
of polynomials in the generators $a^{\alpha}$ and $b^{\alpha}$, where $\alpha=0,1$,
with coefficients in $G$ (see Eq. (\ref{defg})), satisfying the relations%
\bee
\label{RH}     \nonumber
R_k a^\alpha=  - {\lambda }^k b^\alpha R_k ,&\qquad&
R_k b^\alpha = - {\lambda }^{-k} a^\alpha R_k , \\
\label{SH}     \nonumber
S_k a^\alpha =    {\lambda }^{-k} a^\alpha S_k ,&\qquad&
S_k b^\alpha =   {\lambda }^k     b^\alpha S_k ,
\eee
\bee\label{HH}     \nonumber
\left[a^\alpha,\,b^\beta\right]&=&\varepsilon^{\alpha\beta}
\left( 1+n \nu  L_0 \right), \nn
\left[a^\alpha,\,a^\beta\right]&=&\varepsilon^{\alpha\beta}
 n \nu L_1  ,\nn
\left[b^\alpha,\,b^\beta\right]&=&\varepsilon^{\alpha\beta}
n \nu L_{-1},
\eee
where $\varepsilon^{\alpha\beta}$ is the skew-symmetric tensor with
 $\varepsilon^{0\,1}=1$.

Defining the parity on \HG  \ by setting
\be     \nonumber
\varepsilon(a^\alpha)=\varepsilon(b^\alpha)=1,\ \ \ \varepsilon(R_k)=\varepsilon(S_k)=0,
\ee
we turn this algebra into a superalgebra.

Introduce a new parameter of the algebra \HG:
\be\label{3.5}
\mu:=n {\nu} ,
\ee
and rewrite the relations between the generating elements of
 \HG \ and elements of the group algebra $G$:
\bee
\label{LH}     \nonumber
&L_p a^\alpha = -   b^\alpha L_{p+1} ,\qquad
&L_p b^\alpha = -   a^\alpha L_{p-1} ,\\
\label{QH}
&Q_p a^\alpha = a^\alpha Q_{p+1} ,\qquad
&Q_p b^\alpha = b^\alpha Q_{p-1} ,\\
\label{GG}
&L_k L_l =  \delta_{k+l} Q_l, \qquad &L_k Q_l =  \delta_{k-l} L_l,\nn
&Q_k L_l =  \delta_{k+l} L_l, \qquad & Q_k Q_l =  \delta_{k-l} Q_l,
\ \ \mbox{where }\delta_k\defeq \delta_{k0},     \nonumber
\eee
\bee\label{HH2}
&\left[a^\alpha,\,b^\beta\right]=\varepsilon^{\alpha\beta}
\left( 1+ \mu  L_0 \right), &\nn
&\left[a^\alpha,\,a^\beta\right]=\varepsilon^{\alpha\beta}
\mu   L_1 ,&\nn
&\left[b^\alpha,\,b^\beta\right]=\varepsilon^{\alpha\beta}
\mu   L_{-1}.&
     \nonumber
     \eee


\section{Subalgebra of singlets}\label{sl2}

Consider the elements%
\footnote
{Here the brackets $\{\cdot ,\cdot\}$ denote anticommutator.
}
$
T^{\alpha\beta}:=\frac 1 2 (\{a^\alpha,\, b^\beta \}+\{b^\alpha,\, a^\beta \})
$
of the algebra \HG
\ and the inner derivations they generate:
\be     \nonumber
D^{\alpha\beta}:\ \ f\mapsto \left[f,T^{\alpha\beta}\right] \ \
\text{for any } f\in\HH.
\ee
It is easy to verify that the linear span of these derivations
is isomorphic to the Lie algebra $sl_2$.

\definition
A \emph{singlet} is any element $f\in \HH$ such that $[f,\, T^{\alpha\beta}]=0$ for all $\alpha,\beta$.
The subalgebra $H^0 \subset \HH$ consisting of all singlets of the algebra
\HG \ is called  the \emph{subalgebra of singlets}.

One can consider the algebra \HG \ as an $sl_2$-module and decompose it into
the direct sum of irreducible submodules.

Observe, that any $\varkappa$-trace is identically zero on all irreducible
$sl_2$-submodules of \HG \ except singlets.

Let the skew-symmetric tensor $\varepsilon_{\alpha\beta}$ be normalized so that
$\varepsilon_{0\,1}=1$. We set
\bee     \nonumber
\m:=\frac 1 {4i} \sum_{\alpha,\beta=0,1}
\left(\{a^\alpha,\, b^\beta \}-\{b^\alpha,\, a^\beta \}\right)\varepsilon_{\alpha\beta}.
\eee

It is easy to prove the following fact:
\proposition \label{1.3.7}
{\it The subalgebra of singlets $H^0$ of the algebra \HG \ is the algebra of polynomials in
the element $\m$ with coefficients in the group algebra  $\G$.}

In what follows we need the commutation relations of the singlet $\m$ with generators of
the algebra \HG:
\bee     \label{comsin}
[\m ,\,Q_p] = [\m ,\,S_k] = [T^{\alpha\beta},\,\m]=0,     \\
\m L_p =-L_p\m, \qquad \m R_k =-R_k\m,          \nonumber  \\
(\m -i\mu L_0 )a^\alpha = a^\alpha (\m + i+i \mu L_0).       \nonumber
\eee


\theorem\label{i0}
{\it
Let $\mathcal I$ be a proper ideal in the algebra \HG,
$\mathcal I_0 :=\mathcal I \bigcap H^0$.
Then there exist  nonzero polynomials
$\phi_k^0\in \mathbb C[\m]$, where $k=0,...,n-1$, such that
$\mathcal I_0$ is the span over $\mathbb C[\m]$ of the elements
\bee
\phi_k^0 (\m)Q_k, \qquad \ \phi_{n-k}^0 L_k, \text{ where } k=0,\,...\,,\,n-1
\ \text{and }\ \phi_n^0:=\phi_0^0.
\eee
}

Before proving Theorem \ref{i0}, we formulate and prove several propositions.

\proposition\label{ne0}
{\it
If $\mathcal I\subset \HH$ is a proper ideal,
then
$\mathcal I_0 =\mathcal I \bigcap H^0$ is a proper ideal in
$H^0$.
}

\begin{proof}
First,  note that $\mathcal I_0 \ne H^0$ because $\mathcal I$ does not contain unit.

Second, to prove that $\mathcal I_0 \ne 0$, we consider a nonzero element $g\in \mathcal I$.
The $sl_2$-action  on $g$ generates an invariant subspace $\mathcal F\subset \mathcal I$,
which can be decomposed into sum of invariant subspaces,
$\mathcal F=\bigoplus_s \mathcal F^s$,
where $\mathcal F^s \subset \mathcal I$ is a direct sum of irreducible  $sl_2$-modules of
spin $s$ (and dimension $2s+1$).

We further consider the highest-weight vector $f\in \mathcal F^s$ and the set of elements
\{${fQ_p \mid p=0,...,n-1}$\}, belonging to the ideal $\mathcal I$.
Not all these elements are equal to zero because $\sum_p fQ_p  =f$.
Let $fQ_p\ne 0$ and let it be of degree $N$.
We consider the highest-degree part of the polynomial $fQ_p$, which has the form
$f_Q Q_p+f_L L_p$, where $f_Q$ and $f_L$ are homogeneous polynomials in $a^\alpha$, $b^\alpha$
of degree $N$.
We can assume that $f_Q\ne 0$ (otherwise we can take an element $(f_L L_p)L_{-p}=f_L Q_{-p}\ne 0$)
and consider the polynomial ${\m fQ_p + fQ_p \m=2\m f_Q Q_p}$.
The highest-degree terms of this polynomial have the form
\be     \nonumber
\sum_{k=0}^{2s} c_k (a^1)^k (b^1)^{2s-k} \m^{(N+2)/2-s} Q_p,
\ee
where $N+2$ is the degree of the homogeneous polynomial $\m f_Q + f_Q \m$.
Let $c_k\ne 0$. Let us consider an element $(b^0)^k (a^0)^{2s-k} (\m f Q_p + f Q_p \m  )$ from the ideal
$\mathcal I$
and the invariant subspace that it generates under the  $sl_2$-action. It contains
a nonvanishing subspace of singlets.
\end{proof}

\definition
For each $p=0, ... , 2m$, we define the ideals $\mathcal J_p$ and $\mathcal J^p$
in the algebra $\mathbb C[\m]$,
by setting
\be     \nonumber
\mathcal J_p := \{f\in \mathbb C[\m] \mid \ f(\m)Q_p \in \mathcal I    \},
\qquad \mathcal J^p := \{f\in \mathbb C[\m] \mid \ f(\m)L_p \in \mathcal I    \}.
\ee

\proposition\label{8.13}
{\it We have
$\mathcal J_p=\mathcal J^{-p}$.
}

\begin{proof}
It follows from the identities $f(\m)Q_p L_{-p}=f(\m) L_{-p}$ and $f(\m) L_{-p}L_p=f(\m) Q_{p}$.
\end{proof}

\proposition\label{9.4}
{\it We have $\mathcal J_p\ne 0$ for any $p=0,...,2m$.
}

\begin{proof}
Let us consider a nonzero element $f\in \mathcal I_0$.

By Proposition \ref{1.3.7}, $f=\sum_p (\phi_p(\m) Q_p + \psi_p(\m) L_{-p})$.
Obviously, there exists a $p$ such that either $\phi_p\ne 0$ or $\psi_p \ne 0$.
So, at least one of the elements $\m Q_p f+ Q_pf\m =2\m\phi_p(\m) Q_p \in \mathcal I_0$ and
$\m Q_p f- Q_p f\m =2\m\psi_p(\m) L_{-p} \in \mathcal I_0$ is nonzero.
Hence, $\mathcal J_p \ne 0$.

Further, we prove that if
$\mathcal J_p \ne 0$, then $\mathcal J_{p+1} \ne 0$, and therefore  $\mathcal J_k \ne 0$
for $k=0,1,...n-1$.

Let $g\in \mathcal J_p$, $g\ne 0$.
Then $gQ_p\in \mathcal I$, and the element
$\tilde g:= \varepsilon _{\alpha\beta}b^\alpha gQ_p a^\beta\in \mathcal I$ is also nonzero.

By relation (\ref{QH}), $\tilde g = \varepsilon _{\alpha\beta}b^\alpha g a^\beta Q_{p+1}$,
with $\tilde g\in \mathcal I_0$, and hence
$\tilde g = \sum_k(\phi_k(\m) Q_k + \psi_k(\m) L_{-k})$ by Proposition \ref{1.3.7}.
Because $0\ne\m\tilde g Q_{p+1}+\tilde g Q_{p+1}\m \in \mathcal I_0$, as can be verified,
we have ${\m\phi_{p+1}(\m)\ne 0}$, and $\m\phi_{p+1}(\m)Q_{p+1}\in \mathcal I_0$,
i.e., $\m\phi_{p+1}(\m) \in \mathcal J_{p+1}\ne 0$.
\end{proof}

Since $\mathbb C[\m]$ is a principal ideal ring, we have the following statement:

 \corollary \label{col} {\it For any
$p=0,...,2m$, there exists a nonzero polynomial $\phi_p^0\in \mathbb C[\m]$
such that $\mathcal J_p =\phi_p^0 \mathbb C[\m]$.
}

Theorem \ref{i0} evidently follows from Corollary \ref{col}.


\section{Generating functions of $\varkappa$-traces}\label{genfuni2}

For each $\varkappa$-trace $\spur$ on \HG, one can define the following
set of generating functions, which allows one to calculate the $\varkappa$-trace
of arbitrary element in $H^0$ via finding the derivatives with respect to parameter $t$ at
zero:
\bee\label{1.33}
&& F_p^{\spur}(t):=\spur(\exp(t(\m-i\Lo))Q_p),
  \\
&& \Psi_p^{\spur}(t):=\spur(\exp(t\m)L_p),     \nonumber
\eee
where $p=0,...,2m$%
.

Since $L_0 Q_p=0$ for any $p\ne 0$,
it follows from the definition Eq. (\ref{1.33}) that
\bee\nonumber
F_p^{\spur}(t)&=& \spur(\exp(t \m)Q_p) \text{ if } p\ne 0,
\nn
F_0^{\spur}(t)&=&\spur(\exp(t(\m-i\mu L_0))Q_0).
\nonumber
\eee

It is easy to find $\Psi_p^{\spur}$ for $p\ne 0$.
Since $\m L_q =-L_q\m$ for any $q=0,...,2m$,
we have
\be
\Psi_q^{\spur}(t)=\spur(\exp(t\m)L_q)=\spur(L_q).
\ee
Next, since $\spur(R_k)$ does not depend on $k$, we have
$\spur(L_p)=0$ for any $p\ne 0$ and
\be\label{psi0}
\Psi_p^{\spur}(t)\equiv 0 \text{ for any } p \ne 0.
\ee

The value of $\spur(L_0)$ will be calculated later, in Section \ref{GAi2-sp}.

We consider also the functions $\Phi_p^{\spur}(t):=\spur(\exp(t(\m+i\Lo))Q_p)$.
It is easily verified, by expanding the exponential in a series, that
these function are related with the functions $F_p^{\spur}$
by the formula
\be     \nonumber
\Phi_p^{\spur}(t)=F_p^{\spur}(t)+2i\Delta_p^{\spur}(t), \\
\text{where } \Delta_p(t)^{\spur}=\delta_p \sin(\mu t)\spur(L_0).
     \nonumber
\ee

The form of generating functions is related with (non)degeneracy of bilinear form
$B_{\spur}$ by Proposition \ref{1.7.8} below.


\section{Degeneracy conditions for the $\varkappa$-trace}\label{degcond}

\proposition\label{1.7.8}
{\it The $\varkappa$-trace on the algebra \HG \ is degenerate  if and only if the
generating functions $F_p^{\spur}$ defined by  formula (\ref{1.33})
have the following form
\bee\label{wid1}
F_p^{\spur}(t)=\sum_{j=1}^{j_p} \exp(t \omega_{j,p})\varphi_{j,p}(t),
\eee
where $\omega_{j,p}\in \mathbb C$ and $\varphi_{j,p}\in \mathbb C [t]$.
}

\begin{proof}

{\bf Sufficiency.}
Let the
functions $F_p^{\spur}$ defined by Eq. (\ref{1.33})
have the form (\ref{wid1}).

We introduce the polynomials $D_p\in \mathbb C[x]$ by the formulas
\bee     \nonumber
&& D_p(x):=\prod_{j=1}^{j_p}(x- \omega_{j,p})^{1+\deg \varphi_{j,p}} \qquad \text {for } p\ne 0, \\
&& D_0(x):=\prod_{j=1}^{j_0}(x^2- \omega_{j,0}^2)^{1+\deg \varphi_{j,0}}.       \nonumber
\eee
By definition, these polynomials satisfy the conditions $D_p(\frac d {dt} ) F_p^{\spur}(t) =0$ for
any $p$.
Besides, introduce the polynomial $\tilde D_0$
by setting
\bee     \nonumber
\tilde D_0(x^2)=D_0(x).
\eee

Since the $\varkappa$-trace $sp$ we consider is non-zero, there exists
a $p$ such that $F_p^{\spur}\ne 0$.

Now, we see, that if
$F_p^{\spur}\ne 0$ for some $p\ne 0$, then the element $D_p(\m)Q_p \in\HH$ is
a null-vector of the bilinear form $B_{\spur}$;
we also see that
if
$F_0^{\spur}\ne 0$, then the element ${\widehat D(\m)Q_0 := {\m^2\tilde D_0(\m^2-\mu^2)Q_0 \in\HH}}$
is a null-vector of the bilinear form $B_{\spur}$.

Indeed, if $f\in \HH$ belongs to a nonsinglet irreducible $sl_2$-module,
then $\spur(D_p(\m)Q_p f)=0$ for any $p\ne 0$ and $\spur (\widehat D_0(\m)Q_0 f)=0$.
If $f\in H^0$, then
$f=\sum_q (f_q(\m)Q_q + g_q(\m)L_q)$
and, taking in account Eq. (\ref{psi0}),
\[
\spur(D_p(\m)Q_p f)=\spur(D_p(\m) Q_p f_p)=
f_p\left(\frac d {dt}\right) D_p\left(\frac d {dt}\right)F_p^{\spur}|_{t=0}=0
\text{ for } p \ne 0.
\]

Further,
let us decompose the polynomial $f_0$ in the sum of even and odd polynomials:
\be     \nonumber
f_0(\m)=f_0^+(\m^2)+\m f_0^-(\m^2).
\ee
Since $\spur (\m^k Q_0)=0$ when $k$ is odd%
\footnote
{Indeed,
\[
\spur(\m^k Q_0)=\spur(\m^k L_0 L_0)=\spur(L_0\m^k L_0)=\spur((-1)^k \m^k L_0 L_0)=\spur((-1)^k \m^k Q_0).
\].
},
since
\[
\spur(\m^2\tilde D_0(\m^2-\mu^2)Q_0 g_0 L_0)=0
\]
and
\[
\frac {d^2} {dt^2} \exp(t(\m-i\mu L_0))=
\exp(t(\m-i\mu L_0))(\m^2-\mu^2 Q_0),
\]
it follows that
\bee     \nonumber
\!\!\!\!\!\!\!\!\!\!\!\!\!\!\!\!\!\!\!\!\!\!\!\!\!\!\!\!\!
&& \spur(\m^2\tilde D_0(\m^2-\mu^2)Q_0 f)=\spur(\m^2\tilde D_0(\m^2-\mu^2)Q_0 f_0^+(\m^2))=\\
&&
=(\frac {d^2} {dt^2}+\mu^2)\,f_0^+ \! (\frac {d^2} {dt^2} +\mu^2)
\,\tilde D_0 \! (\frac {d^2} {dt^2})\,F_0(t)|_{t=0}=\nn
&&
=(\frac {d^2} {dt^2}+\mu^2)\,f_0^+ \! (\frac {d^2} {dt^2} +\mu^2)
\, D_0 \! (\frac {d} {dt})\,F_0(t)|_{t=0}=
0.
     \nonumber
\eee

Thus, the sufficiency of Proposition \ref{1.7.8} is proved.

{\bf Necessity.}
We now prove  that if the $\varkappa$-trace is degenerate, then
there exist polynomials $D_p\in \mathbb C[x]$
such that $D_p(\frac d {dt})F_p(t)=0$ for $p=0,...,2m$,
and therefore the generating functions $F_p$ have the form (\ref{wid1}).

Let an ideal $\mathcal I \subset \HH$ consist of null-vectors
of the bilinear form $B_{\spur}$. Then $\mathcal I_0$ consists of
singlet null-vectors, and the vectors $\phi_k^0(\m)Q_k$ and $\phi_k^0(\m)L_{2m+1-k}$
defined by the conditions of Theorem \ref{i0} generate an ideal $\mathcal I_0$ in $\HH_0$.

Let $p\ne 0$.
Then
\be     \nonumber
0\equiv \spur(\phi_p^0(\m) Q_p e^{t\m}Q_p)=\phi_p^0 \! \left(\frac d {dt}\right)F_p(t)
\ee
and, therefore, the function $F_p$ has the form (\ref{wid1}).

Further, we consider the null-vector $\phi(\m^2)Q_0$ of the bilinear form $B_{\spur}$,
where ${\phi(\m^2):= \phi_0^0(\m) \phi_0^0(-\m)}$.
We note that
\be     \nonumber
\frac {d^2}{dt^2}F_0=
\spur \left( e^{t(\m-i\mu L_0)}(\m-i\mu L_0)^2 Q_0 \right)=
\spur \left( e^{t(\m-i\mu L_0)}(\m^2-\mu^2) Q_0 \right),
\ee
hence
\be     \nonumber
\spur \left( e^{t(\m-i\mu L_0)}\m^2 Q_0 \right)=
\left(\frac {d^2}{dt^2}+\mu^2\right)F_0
\ee
and
\be     \nonumber
0\equiv \spur \left(e^{t(\m-i\mu L_0)}\phi_0^0(\m) \phi_0^0(-\m)Q_0 \right)=
\spur \left( e^{t(\m-i\mu L_0)}\phi(\m^2)Q_0\right)=
\phi\!\left(\frac {d^2}{dt^2}+\mu^2\right)F_0(t),
\ee
i.e., the function $F_0$ also has the form (\ref{wid1}).
\end{proof}


\section{Equations for the generating functions $F_p^{\spur}$}\label{genfun}

Let us differentiate the generating function $F_p^{\spur}$:
\be     \nonumber
\frac d {dt} F_p^{\spur}(t)=\spur\left(e^{t(\m-i\Lo)}(\m-i\Lo)Q_p\right)=
\spur\left(e^{t(\m-i\Lo)}(-ia^\alpha\varepsilon_{\alpha\beta}b^\beta+i)Q_p\right).
\ee
The second equality here holds because
\be     \nonumber
\m=-ia^\alpha\varepsilon_{\alpha\beta}b^\beta+i(1+\Lo).
\ee
Next,
\bee
&& \spur\left(e^{t(\m-i\Lo)}(-ia^\alpha\varepsilon_{\alpha\beta}b^\beta)Q_p\right)=
\spur\left(a^\alpha e^{t(\m+i\Lo)}(-i\varepsilon_{\alpha\beta}b^\beta)Q_p\right)=  \nn
&& \,= \varkappa\spur\left( e^{t(\m+i+i\Lo)}(-i\varepsilon_{\alpha\beta}b^\beta a^\alpha) Q_{p+1}\right)=
\varkappa\spur\left( e^{t(\m+i+i\Lo)}(\m + i +i\Lo) Q_{p+1}\right)=    \nn
&&\, =\varkappa\frac d {dt} \left(e^{it}\Phi_{p+1}(t)\right).
     \nonumber
\eee
Thus, we obtain a system of differential equations for the generating functions:
\bee\label{eqgenfun}
\frac d {dt} F_p^{\spur} - \varkappa e^{it}\frac d {dt}F_{p+1}^{\spur}=iF_p^{\spur} +
\varkappa ie^{it}F_{p+1}^{\spur}+2 \varkappa i\frac d {dt}\left(e^{it}\Delta_{p+1}^{\spur} \right).
\eee

The initial conditions for this system are:
 \[
F_p^{\spur}(0) = \spur(Q_p).
\]

To solve the system (\ref{eqgenfun}), we consider its Fourier transform.
Let
\bee     \nonumber
&& \lambda:=e^{2\pi i/(2m+1)}, \\
&& G_k^{\spur}:= \sum_{p=0}^{2m} \lambda^{kp}F_p^{\spur}, \text{ where } k=0,...,2m,
\label{fur}
\\
&& \widetilde{\Delta}_k^{\spur}:=\sum_{p=0}^{2m} \lambda^{kp}\Delta_{p+1}^{\spur}=
\lambda^{-k}\left( \sin(\mu t)\spur(L_0)    \right),
 \text{ where } k=0,...,2m.
     \nonumber
\eee
For the functions $G_k^{\spur}$, we then obtain the equations
\be\label{Geq}
\frac d {dt} G_k^{\spur} = i\frac {\lambda^k+ \varkappa e^{it}}{\lambda^k- \varkappa e^{it}} G_k^{\spur} +
\frac {2i \varkappa \lambda^k }{\lambda^k- \varkappa e^{it}} \frac d {dt}\left(e^{it}\widetilde{\Delta}_k^{\spur}\right)
\ee
with the initial conditions
\be     \label{init}
G_k^{\spur}(0)=\spur(S_k).
\ee

We choose the solution of the system (\ref{Geq}) in the form:
\bee\label{GK}
G_k^{\spur}(t)= \frac {\varkappa e^{it}} {(\varkappa e^{it}-\lambda^k)^2} \lambda^k g_k^{\spur}(t),
\eee
where
\bee\label{1.47g}
g_k^{\spur}(t)=\left(\frac{2}{\mu}(\cos (t\mu) -1)
+{2i}\lambda ^{-k}(\lambda ^{k}-\varkappa e^{it})\sin (t\mu)
\right)\spur(L_0)
+\varkappa\lambda ^{-k}(\varkappa -\lambda ^{k})^{2} \spur(S_k).
\eee

Evidently, this solution satisfies initial condition (\ref{init}) for each $\varkappa$
and $k$, except the case $\varkappa=+1$ and $k=0$.


If $\varkappa = +1$ and $k=0$, then the expression Eq. (\ref{GK})
for $G_0^{tr}$ has a removable singularity at $t=0$.
In this case, we consider the condition $\lim_{t \rightarrow 0} G_0^{tr}(t) = tr(S_0)$
instead of $G_0^{tr}(0)=tr(S_0)$.

When $\varkappa=+1$ the solution (\ref{GK}) -- (\ref{1.47g}) gives
\be
G_0^{tr}(t)= \frac { e^{it}} {( e^{it}-1)^2}
\left(\frac{2}{\mu}(\cos (t\mu) -1)
+{2i}(1- e^{it})\sin (t\mu)
\right)tr(L_0)
\ee
and one can easily see that
\be
\lim_{t \rightarrow 0} G_0^{tr}(t) = -\mu tr(L_0).
\ee

It id shown in Subsection \ref{GAi2-tr}
using Ground Level Conditions, that if $\varkappa=+1$,
then
\be
tr(S_0) = -\mu tr(L_0)
\ee
for any trace $tr$ on \HG.

So, the $G_0^{tr}(t)$ satisfies the initial conditions (\ref{init}) also.


For the case $ \varkappa = -1$, the $\varkappa$-trace is a supertrace (see \cite{KT}).
In this case, the $m+1$ values $str(S_k)=str(S_{2m+1-k})$ for $k=0,...,m$ completely define
the supertrace on \HG \ (see \cite{stek}).

For the case $ \varkappa = +1$, the $\varkappa$-trace is a trace (see \cite{KT}).
In this case, the $m$ values ${tr(S_k)=tr(S_{2m+1-k})}$ for $k=1,...,m$ completely define
the trace on \HG \ (see \cite{stek}).
The value $tr(S_0)$ linearly depends on parameters $tr(S_k)$, where $k=1,...,m$,
and it is found in Subsection \ref{GAi2-tr} (see Eqs. (\ref{s0}) -- (\ref{X})).

\section{Values of the $\varkappa$-trace on $\G$}\label{GAi2-sp}

To use the generating functions (\ref{GK}), we need
to express the values $\spur (S_k)$ and $\spur(L_0)$
via some independent parameters which completely define a the $\varkappa$-trace.

The results are different for traces ($\varkappa = +1$) and
for supertraces ($\varkappa = -1$).
First, we express $sp(L_0)$ via $sp(S_k)=sp(S_{2m+1-k})$, where $k=1\,,...\,,\,m$
if $\varkappa=+1$ and $k=0, \,1\,,...\,,\,m$
if $\varkappa=-1$.

Let
\be
c_k^\alpha:= a^\alpha - \varkappa \lambda^k b, \text{ so }\ R_k c_k^\alpha=\varkappa c_k^\alpha R_k.
\ee

We consider the chain of equalities
\be
\spur(c_k^0 c_k^1 R_k) = \varkappa \spur( c_k^1 R_kc_k^0)=\varkappa^2\spur(c_k^1 c_k^0 R_k),
\ee
which results in
\be\label{GLE1}
\spur([c_k^0,\, c_k^1] R_k)=0.
\ee
The conditions like (\ref{GLE1}) are called \emph{ Ground Level Conditions}
in \cite{KV}, \cite{KT}.
It follows from (\ref{GLE1}) that
\be\nonumber
-2\lambda^k \varkappa
\spur\left(R_k-\frac \mu 2 \varkappa (\lambda^{-k} L_1 - 2\varkappa L_0 +\lambda^{k} L_{-1})R_k)
 \right)=0,
\ee
which gives
\be\label{rk}
\spur (R_k) = - \frac {2\mu} {2m+1}
\left(\frac
{1+\varkappa} 2 X^{\spur}
+
\frac {1-\varkappa} 2 Y^{\spur}
\right),
\ee
where
\bee
&& X^{\spur}:= \sum_{r=1}^{2m}
\sin^2\left(\frac {\pi r}{2m+1}\right) \spur (S_r),\\
&& Y^{\spur}:= \sum_{r=0}^{2m}
\cos^2\left(\frac {\pi r}{2m+1}\right) \spur (S_r).
\eee

Below we consider these values for the traces and supertraces separately.

\subsection{Values of the traces on $\G$}\label{GAi2-tr}

The group $I_2(2m+1)$ has $m$ conjugacy
classes without the eigenvalue +1 in the spectrum:
\be     \nonumber
\{S_{p},S_{2m+1-p}\}, \text { where } p=1,...,m.
\ee

By Theorem 2.3 in \cite{KT}, the values of the trace on these conjugacy classes
\be\label{1.2.3}     \nonumber
s_k := tr(S_k),\ \ \text{ where } s_{2m+1-k}=s_k, \ \ k=1,...,m,
\ee
are arbitrary and
completely defines the trace on the algebra \HG,
and therefore the dimension of the space of traces is $m$.

Further,
the group $I_2(2m+1)$ has one conjugacy class with one eigenvalue +1 in its spectrum:
\be     \nonumber
\{R_{1},\,...\,,\, R_{2m+1}\}.
\ee
The value of $tr(R_k)$ is expressed via $s_k$ by formula (\ref{rk}).

Besides, the group $I_2(2m+1)$ has
one conjugacy class with two eigenvalues +1 in its spectrum: $\{S_0\}$.

The traces on conjugacy classes with two eigenvalues +1 in the spectrum also
can be calculated using Ground Level Conditions (see \cite{KT}):
\bee
tr([a^0,\, b^1]S_0)=0,      \nonumber
\eee
which gives
\bee  \label{R0}
\label{s0}
tr(S_{0}) = 2 \nu^{2}(2m+1) X^{tr},
\eee
where
\bee\label{X}
X^{tr}:= \sum_{l=1}^{2m}s_{l}\sin ^{2}\left(\frac{2\pi l}{2m+1}\right).
\eee

We also note that
\bee     \nonumber
tr(L_0)=-\frac {2\mu} {2m+1} X^{tr},\quad
tr (L_p)=0 \text{ for } p\ne 0, \quad
tr(S_0)=-\mu tr(L_0).
     \nonumber
\eee


\subsection{Values of the supertraces on $\G$}\label{GAi2-str}

The group $I_2(2m+1)$ has $m+1$ conjugacy
classes without the eigenvalue -1 in the spectrum:
\be     \nonumber
\{S_{0}\},\
\{S_{p},S_{2m+1-p}\}, \text { where } p=1,...,m.
\ee
By Theorem 2.3 in \cite{KT}, the values of the supertrace on these conjugacy classes
\be\label{1.2.3-}     \nonumber
u_k := str(S_k)=str(S_{2m+1-k}),\ \ \text{ where }  \ \ k=0,...,m,
\ee
are arbitrary parameters that
completely define the supertrace $str$ on the algebra \HG,
and therefore the dimension of the space of supertraces is $m+1$.

Besides,
the group $I_2(2m+1)$ has one conjugacy class with one eigenvalue -1 in the spectrum:
\be     \nonumber
\{R_{1},\,...\,,\, R_{2m+1}\}.
\ee

The supertraces of the conjugacy class with eigenvalue $-1$ in  its spectrum
is calculated via Ground Level Conditions in Section \ref{GAi2-sp}.
These conditions give
\bee     \nonumber
str(R_{k}) = -2\nu Y^{str},
&\ \ & k=0,1,...,2m,     \nonumber
\eee
where
\bee     \nonumber
Y^{str}:= \sum_{r=0}^{2m} u_{r}\cos ^{2}
\left( \frac{\pi r}{2m+1} \right).
\eee


\section{Singular values of the parameter $\mu$}\label{spec}

We now find  the values of the parameter $\mu$ for which there exists
a nonzero $\varkappa$-trace $\spur$, i.e., the values $\spur(S_k)$ such that the
the generating functions $F_p$ (\ref{1.33}) have the form (\ref{wid1}).
Since the functions $G_k$ (\ref{fur}) are linear combinations of the functions $F_p$,
and vice versa, the algebra \HG \ has a degenerate  $\varkappa$-trace if and only if
the functions $G_k$ (\ref{fur}) have the form (\ref{wid1}) also.

In particular, it is necessary that
the numerator of the expression (\ref{GK}) contains all the zeros of
the denominator of the expression.

The denominator of the function $G_k$ is equal to
\be\label{znam}     \nonumber
(e^{it}-\varkappa\lambda^k)^2
\ee
and have doubled zeros at
\be     \nonumber
t_{k,l}=\frac {2\pi} n k+ 2\pi l +\pi \theta, \text{ where }
\quad l=0,\, \pm 1,\, \pm 2, \, ...
\ee
and
\be
\theta = \left\{
\begin{array}{l}
0 \text{ if } \varkappa=1 \\
1 \text{ if } \varkappa=-1.
\end{array}
\right.
\ee

It is easy to check that $\frac {d}{dt} g_k^{\spur} (t_{k,l_k})=0$ for each $k=0,...,2m$ and each integer $l_k$.

The equalities $g_k^{\spur} (t_{k,l_k})=0$ can be considered as a system of linear equations for
the values $tr(S_k)=tr(S_{n-k})$, where $k=1,...,2m$ if $\varkappa=1$, and for the values
$str(S_k)=str(S_{n-k})$, where $k=0,...\,,\,n$ if $\varkappa=-1$:
\bee\label{147}
g_k^{\spur}(t_{k,\,l_k})=\frac{2}
{\mu} \left(\cos (t_{k,\,l_k}\mu) -1 \right)
\spur(L_0)
+\varkappa\lambda ^{-k}(\varkappa -\lambda ^{k})^{2} \spur(S_k)=0.
\eee
Our goal is to find the $\mu$, such that the system (\ref{147}) has nonzero solutions.

Note that $\spur(L_0)\ne 0$ otherwise the $\varkappa$-trace would be zero.
We consider the subsystem of two equations with $l_k=0$%
:
\bee\label{2eq}
\frac{2}
{\mu}(\cos ((\frac {2\pi k} n  +\pi \theta)\mu) -1)
\spur(L_0)
+\varkappa\lambda ^{-k}(\varkappa -\lambda ^{k})^{2} \spur(S_k)=0,
\\
\label{2eq-}
\frac{2}
{\mu}(\cos ((\frac {2\pi(n-k)} n  +\pi \theta)\mu) -1)
\spur(L_0)
+\varkappa\lambda ^{k-n}(\varkappa -\lambda ^{n-k})^{2} \spur(S_{n-k})=0.
\eee
Since $\frac {\varkappa(\varkappa -\lambda ^{k})^{2}} {\lambda ^{k}}=
\frac {\varkappa(\varkappa -\lambda ^{n-k})^{2}} {\lambda ^{n-k}}
$ and $\spur(S_{n-k})=\spur(S_{k})$,
Eqs. (\ref{2eq}) -- (\ref{2eq-}) imply that
\be
\cos ((\frac {2\pi k} n +\pi \theta)\mu)-\cos ((\frac {2\pi(n-k)} n  +\pi \theta)\mu)=0
\ee
or
\be\label{sol1}
\sin(\pi\mu(1+\theta))
\sin(\frac {2k-n}{n}\pi\mu)
=0.
\ee
Eq. (\ref{sol1}) implies that
\be
\mu=\frac {z}{1+\theta}, \text{ where } z\in \mathbb Z.
\ee

Next, we consider the two cases separately:

A) $\mu\in \mathbb Z$, $\varkappa = \pm 1$,

B) $\mu = z+\frac 1 2$, where $z\in \mathbb Z$, $\varkappa = - 1$.

To solve the case A), we note that Eq. (\ref{147}) gives for $\mu$ integer:
\be\label{sum}
0=\sum _{k=0}^{n-1} g_k^{\spur}(t_{k,\, l_k})=\frac 2 \mu
\sum _{k=0}^{n-1}\cos(\frac {2k\pi}{n}\mu)(-1)^{\theta\mu} \spur(L_0).
\ee
Since $\spur(L_0)\ne 0$, Eq. (\ref{sum}) gives the following restriction on the integer $\mu$:
\be\label{suum}
\sum _{k=0}^{n-1}\cos(\frac {2k\pi}{n}\mu)=0,
\ee
i.e.,
\be
\mu \in \mathbb Z \setminus n\mathbb Z.
\ee

Now consider the case B), i.e., $\varkappa =-1$, $\theta=1$, $\mu=z+ \frac 1 2$, where $z\in\mathbb Z$,
namely, consider the following two equations of the system (\ref{147}):
\bee
g_k^{str}(t_{k,\,0}) & = & \frac{2}
{\mu}(\cos (
\frac {2\pi k z} n+
\frac {\pi k}n + \pi z +\frac \pi 2
) -1)
\spur(L_0)
-\frac {(1+ \lambda^k )^{2}} {\lambda^k}  str(S_k)=0,\nn
g_k^{str}(t_{k,\,1}) & = &
\frac{2}
{\mu}(\cos (
\frac {2\pi k z} n+
\frac {\pi k}n + \pi z +\frac \pi 2
+\pi) -1)
\spur(L_0)
-\frac {(1+ \lambda^k )^{2}}{\lambda^k}  str(S_k)=0,
\nonumber
\eee
which give
\be
\cos (
\frac {2\pi k z} n+
\frac {\pi k}n + \pi z +\frac \pi 2
)=0
\ee
or
\be
2z+1=n r \ \text{ for some odd $r$, or $\mu=\frac {nr}{2}$. }
\ee

One easily checks that for every $\mu$ found,
the system (\ref{147}) does not depend on $l_k$ and so has a nonzero solution.

Thus, we have proved the following theorem:

\vskip 2mm
\theorem\label{th2}
{\it
Let $ m\in \mathbb Z$, $m\geqslant 1$
and $n=2m+1$.
Then

1) The associative algebra $H_{1,\nu}(I_2(n))$
has a one-parameter set of nonzero  traces $tr_z$ such that the symmetric invariant bilinear form
$B_{tr_z}(x,y)=tr(xy)$  is degenerate if and only if
$\nu= \frac z n$, where $z\in \mathbb Z \setminus n\mathbb Z$.
These traces are completely defined by their values on $S_k$ for
$k=1,\dots , m$:
\be\label{913}
tr_z(S_k)= \frac \tau {n \sin^{2}(\frac {\pi k}n)}
(1 - \cos(\frac {2\pi k z}{n})), \text{ where } \tau\in \mathbb C,\ \tau\ne 0.
\ee

2) The associative superalgebra $H_{1,\nu}(I_2(n))$
has a one-parameter set of  nonzero supertraces $str_z$ such that the supersymmetric invariant bilinear form
$B_{str_z}(x,y)=str(xy)$  is degenerate if
$\nu= \frac z n$, where $z\in \mathbb Z \setminus n\mathbb Z$.
These supertraces are completely defined by their values on $S_k$ for
$k=0,\dots , m$:
\bee\label{914}
str_z(S_k)= \frac \tau {n \cos^{2}(\frac {\pi k}n)}
(1 - (-1)^{z}\cos(\frac {2\pi k z}{n})), \text{ where } \tau\in \mathbb C,\ \tau\ne 0.
\eee

3) The associative superalgebra $H_{1,\nu}(I_2(n))$
has a one-parameter set of nonzero supertraces $str_{1/2}$ such that the supersymmetric invariant bilinear form
$B_{str_{1/2}}(x,y)=str_{1/2}(xy)$  is degenerate if
$\nu= z + \frac 1 2$, where $z\in \mathbb Z$.
These supertraces are completely defined by their values on $S_k$ for
$k=0,\dots , m$:
\be
str_{1/2}(S_k)= \frac \tau {n \cos^{2}(\frac {\pi k}n)}, \text{ where } \tau\in \mathbb C, \ \tau\ne 0.
\ee

4) For all other values of $\nu$, all nonzero traces and supertraces are nondegenerate.
}

\remark
Theorem \ref{th2} implies that if $z\in \mathbb Z \setminus n\mathbb Z$, then
the trace (\ref{913}) generates the ideal $\mathcal I_{tr_z}$ consisting of null-vectors
of the degenerate form $B_{tr_z}(x,y)=tr_z(xy)$, and simultaneously
the supertrace (\ref{914}) generates the ideal $\mathcal I_{str_z}$ consisting of null-vectors
of the degenerate form $B_{str_z}(x,y)=str_z(xy)$.
A question arises: is it true that $\mathcal I_{tr_z}=\mathcal I_{str_z}$?

\conjecture
$\mathcal I_{tr_z}=\mathcal I_{str_z}$.

Our observation, that the set of coefficients $\omega_{j,p}$
in Eq. (\ref{wid1}) for $F_p^{tr_z}$ is the same as for $F_p^{str_z}$,
is an argument in favor of this conjecture.


\vskip 15mm

\section*{Acknowledgments}
The authors (S.K. and I.T.) are grateful to Russian Fund for Basic Research
(grant No.~${\text{14-02-01171}}$)
for partial support of this work.


\setcounter{equation}{0} \def\theequation{A
\arabic{equation}}

\newcounter{appen}
\newcommand{\appen}[1]{\par\refstepcounter{appen}
{\par\bigskip\noindent\large\bf Appendix. 
\medskip }{\bf \large{#1}}}

\renewcommand{\subsection}[1]{\refstepcounter{subsection}
\vskip 3mm{\centerline{\bf 
A\arabic{subsection}. \ #1}}\vskip 3mm}
\renewcommand\thesubsection{A\theappen.\arabic{subsection}}
\makeatletter \@addtoreset{subsection}{appen}

\renewcommand{\subsubsection}{\par\refstepcounter{subsubsection}
{\bf A\arabic{appen}.\arabic{subsection}.\arabic{subsubsection}. }}
\renewcommand\thesubsubsection{A\theappen.\arabic{subsection}.\arabic{subsubsection}}
\makeatletter \@addtoreset{subsubsection}{subsection}


\newcounter{theor}
\renewcommand{\theorem}{\par\refstepcounter{theor}
           {\bf Theorem A.%
           \arabic{theor}. }}
\renewcommand\thetheorem{A.\arabic{theor}}
\makeatletter \@addtoreset{theor}{section}

\renewcommand{\remark}{\par\refstepcounter{theor}
           {\bf Remark A.%
           \arabic{theor}. }}
\renewcommand\theremark{\thesection.%
     \arabic{theor}}
\makeatletter \@addtoreset{remark}{section}

\renewcommand{\definition}{\par\refstepcounter{theor}
           {\bf Definition A.%
           \arabic{theor}. }}
\renewcommand\thedefinition{\thesection.%
    \arabic{theor}}
\makeatletter \@addtoreset{definition}{section}

\renewcommand{\proposition}{\par\refstepcounter{theor}
           {\bf Proposition 
           A.\arabic{theor}. }}
\renewcommand\theproposition{
  A.\arabic{theor}}
\makeatletter \@addtoreset{proposition}{section}

\renewcommand{\theorem}{\par\refstepcounter{theor}
           {\bf Theorem 
           A.\arabic{theor}. }}
\renewcommand\thetheorem{
   A.\arabic{theor}}
\makeatletter \@addtoreset{theorem}{section}


\appen{The case $H_{1,\nu_1,\nu_2}(I_2(n))$ with $n$ even}\label{App}

Here we, following \cite{tmf}, briefly describe the degenerate traces generating the ideals in
the Symplectic Reflection Algebra $H_{1,\nu_1,\nu_2}(I_2(2m))$.

This algebra has two complex parameters; for every value of these parameters
the algebra has an $m$-dimensional space of traces and, due to presence of the Klein operator,
the isomorphic space of supertraces.

\subsection{The group $I_2(2m)$}

\definition
{The group $I_2(2m)$ is a finite subgroup of $O(2,\mathbb R)$, generated by
the root system $I_2(2m)$.
It consists of
$2m$ reflections $R_k$, acting on $z \in{\mathbb C}$
as follows
\bee\label{ARk}
&& R_k z   = -z^* v_k^2 R_k, \quad k=0,\,...,\, 2m-1
\eee
and $2m$ rotations $S_k:=R_k R_0$, where $S_0$ is the unit in $I_2(2m)$
and $S_m$ is the Klein operator.
As we see from (\ref{ARk}), these elements satisfy the relations
\be\label{ARR}     \nonumber
R_k R_l = S_{k-l},\qquad
S_k S_l = S_{k+l},\qquad
R_k S_l = R_{k-l},\qquad
S_k R_l = R_{k+l}.
\ee
}

Evidently, the $R_{2k}$ belong to one conjugacy class and the $R_{2k+1}$ belong to another class.
The rotations $S_k$ and $S_l$ constitute a conjugacy class if $k+l=2m$.

\definition
\bee\label{ALQbas}
&& L_p:=\frac 1 n \sum_{k=0}^{2m-1} \lambda^{kp} R_k,\qquad
Q_p:=\frac 1 n \sum_{k=0}^{2m-1} \lambda^{-kp} S_k,\\
\label{Alambda}
&& \mbox{ where }\lambda=\exp\left(\frac {\pi i} m \right).
\nonumber
\eee

\subsection{Symplectic reflection algebra
$H_{1,\nu_0, \nu_1}(I_2(2m))$}\label{Asrai2}

\definition
{The symplectic reflection algebra
$\mathcal H :=H_{1,\nu_0, \nu_1}(I_2(2m))$} is an associative algebra
of polynomials in $a^{\alpha},b^{\alpha}$, where $\alpha=0,1$,
with coefficients in $\mathbb C [I_2(2m)]$, satisfying the relations%
\bee
\label{ARH}     \nonumber
R_k a^\alpha=  - {\lambda }^k b^\alpha R_k ,&\qquad&
R_k b^\alpha = - {\lambda }^{-k} a^\alpha R_k , \\
\label{ASH}     \nonumber
S_k a^\alpha =    {\lambda }^{-k} a^\alpha S_k ,&\qquad&
S_k b^\alpha =   {\lambda }^k     b^\alpha S_k ,
\eee
\bee
\label{ALH}     \nonumber
&L_p a^\alpha = -   b^\alpha L_{p+1} ,\qquad
&L_p b^\alpha = -   a^\alpha L_{p-1} ,\\
\label{AQH}
&Q_p a^\alpha = a^\alpha Q_{p+1} ,\qquad
&Q_p b^\alpha = b^\alpha Q_{p-1} ,\\
\label{AGG}
&L_k L_l =  \delta_{k+l} Q_l, \qquad &L_k Q_l =  \delta_{k-l} L_l,\nn
&Q_k L_l =  \delta_{k+l} L_l, \qquad & Q_k Q_l =  \delta_{k-l} Q_l,
\ \ \mbox{where }\delta_k\defeq \delta_{k0},     \nonumber
\eee
\bee\label{AHH2}
&\left[a^\alpha,\,b^\beta\right]=\varepsilon^{\alpha\beta}
\left( 1+ \mu_0  L_0 +  \mu_1   L_{m}\right), &\nn
&\left[a^\alpha,\,a^\beta\right]=\varepsilon^{\alpha\beta}
\left(  \mu_0   L_1 +  \mu_1  L_{m+1} \right),&\nn
&\left[b^\alpha,\,b^\beta\right]=\varepsilon^{\alpha\beta}
\left(  \mu_0   L_{-1} +  \mu_1  L_{m-1} \right),&
     \nonumber
     \eee
where $\varepsilon^{\alpha\beta}$ is the skew-symmetric tensor with $\varepsilon^{0\,1}=1$
and
\be\label{A3.5}
\mu_0:=m ( {\nu_0+\nu_1}) ,\ \ \ \mu_1 :=m ( {\nu_0-\nu_1}).
\ee

The basis elements of Lie algebra  $sl_2$ of inner derivations
$
T^{\alpha\beta}:=\frac 1 2 (\{a^\alpha,\, b^\beta \}+\{b^\alpha,\, a^\beta \})
$
act on \HG \ as follows
\be     \nonumber
 f\mapsto \left[f,T^{\alpha\beta}\right] \ \
\text{for each } f\in\HH.
\ee

Let the skew-symmetric tensor $\varepsilon_{\alpha\beta}$ be such that
$\varepsilon_{0\,1}=1$. Set
\bee     \nonumber
\m:=\sum_{\alpha,\beta=0,1}
\frac 1 {4i} (\{a^\alpha,\, b^\beta \}-\{b^\alpha,\, a^\beta \})\varepsilon_{\alpha\beta}.
\eee

Then
\bee     \nonumber
[\m ,\,Q_p] = [\m ,\,S_k] = [T^{\alpha\beta},\,\m]=0,     \\
\m L_p =-L_p\m, \qquad \m R_k =-R_k\m,          \nonumber  \\
(\m -i(\mu_0 L_0 + \mu_1 L_m) )a^\alpha = a^\alpha (\m + i+i (\mu_0 L_0 + \mu_1 L_m)).       \nonumber
\eee


\subsection{The values of the trace on $\mathbb C[I_2(2m)]$}\label{AGAi2}

The group $I_2(2m)$ has $m$ conjugacy classes without the eigenvalue +1 in their spectra:
\be     \nonumber
\{S_{p},S_{n-p}\}, \text { where } p=1,...,m-1,\;\text{ and also }\{S_{m}\}.
\ee
Due to Theorem 2.3 in \cite{KT}, the values of the trace on these conjugacy classes
\be\label{A1.2.3}
s_k := tr(S_k),\ \ \text{ where } s_{2m-k}=s_k, \ \ k=1,...,m,
\ee
completely define the trace on \HG, and therefore the dimension of the space of traces
is  equal to $m$.

The group $I_2(2m)$ has two conjugacy classes each having one eigenvalue +1 in its spectrum:
\be     \nonumber
\{R_{2l} \mid l=0,...,m-1\},\;\{R_{2l+1} \mid  l=0,...,m-1\},
\ee
and one conjugacy class with two eigenvalues +1 in its spectrum: $\{S_0\}$.

The traces on these conjugacy classes are calculated via Ground Level Conditions
\cite{KT}%
:
\bee
&& tr([c_k^0, \, c_k^1]R_k)=0, \text{ where } c_k^\alpha := a^\alpha - \lambda^k b^\alpha
\text{ are eigenvectors of } R_k,
\quad R_k c_k^\alpha = c_k^\alpha R_k ,\nn
&& tr([a^0,\, b^1]S_0)=0      \nonumber
\eee
and are equal to
\bee     \nonumber
tr(R_{2l})&=&-2\nu _{2}X_1 -2\nu _{1} X_2,\\
tr(R_{2l+1})&=&-2\nu _{1} X_1 - 2\nu _{2} X_2 ,
\\
&\ \ & l=0,1,...,m-1,     \nonumber\\
\label{As0}
tr(S_{0}) &=& 2(\nu _{1}^{2}+\nu _{2}^{2})m
X_1+4\nu _{1}\nu_{2} m X_2,
\eee
where
\bee     \nonumber
&& X_1:= \sum_{l=1}^{m-1}s_{2l}\sin ^{2}\left(\frac{\pi
l}{m}\right),
\\
&& X_2:=\sum_{l=0}^{m-1}s_{2l+1}\sin ^{2}%
\left(\frac{\pi (2l+1)}{2m}\right).     \nonumber
\eee

We note also that
\bee     \nonumber
&& tr(L_0)=-\frac {\mu_0} m (X_1+X_2),\quad tr(L_m)=-\frac {\mu_1} m (X_1-X_2),\quad
tr (L_p)=0 \text{ for } p\ne0,\,m,
\\
&& tr(S_0)=-\mu_0 tr(L_0)-\mu_1 tr(L_m).
     \nonumber
\eee


\subsection{Generating functions of the trace}\label{Agenfuni2}

Set $\mathcal L :=\mu_0 L_0 + \mu_1 L_m$.

For each trace $tr$, we define the following set of generating functions on \HG:
\bee\label{A1.33}
&& F_p(t):=tr(\exp(t(\m-i \mathcal L ))Q_p),
  \\
&& \Psi_p(t):=tr(\exp(t\m)L_p),     \nonumber
\eee
where $p=0,...,2m-1$.
From $\m L_p = - L_p \m$ and definition of the trace it follows  that
\be     \nonumber
\Psi_p(t)=\Psi_p(0).
\ee
%
%
We also consider the functions $\Phi_p(t):=tr(\exp(t(\m+i \mathcal L))Q_p)$
related with the functions $F_p$
by the formula
\be     \nonumber
\Phi_p(t)=F_p(t)+2i\Delta_p(t), \\
\text{where } \Delta_p(t)=\delta_p \sin(\mu_0 t)tr(L_0)+\delta_{m-p} \sin(\mu_1 t)tr(L_m).
     \nonumber
\ee

Analogously to our previous consideration,
one can get the following system of equations
\bee\label{Aeqgenfun}
\frac d {dt} F_p - e^{it}\frac d {dt}F_{p+1}=iF_p +
ie^{it}F_{p+1}+2i\frac d {dt}\left(e^{it}\Delta_{p+1} \right).
\eee


Next, we consider the Fourier transform of (\ref{Aeqgenfun}),
namely, we consider
\bee     \nonumber
&& G_k:= \sum_{p=0}^{2m-1} \lambda^{kp}F_p, \text{ where } k=0,...,2m-1,
\label{Afur}
\\
&& \widetilde{\Delta}_k:=\sum_{p=0}^{2m-1} \lambda^{kp}\Delta_{p+1}=
\lambda^{-k}\left( \sin(\mu_0 t)tr(L_0) + \lambda^{km}\sin(\mu_1 t)tr(L_m)    \right),
 \nn
&& \qquad\qquad\qquad \text{ where } k=0,...,2m-1 \qquad
 \text{ and $ \lambda:=e^{i\pi/m}$,}
     \nonumber
\eee
and obtain the system of equation
\be     \nonumber
\frac d {dt} G_k = i\frac {\lambda^k+e^{it}}{\lambda^k-e^{it}} G_k +
\frac {2i}{\lambda^k-e^{it}} \frac d {dt}\left(e^{it}\widetilde{\Delta}_k\right)
\ee
with initial conditions
\be     \label{Ainit}
G_k(0)=s_k, \text{ where } k=0,...,2m-1,
\ee
and where the $s_k$ are defined by Eq. (\ref{A1.2.3}) for $k=1,\,...,\,2m-1$  and
$s_0 := tr(S_0)$ is defined by Eq. (\ref{As0}).
The value $s_0$ depends linearly on
$s_k$, where $k=1,...,m$ (see Eq. (\ref{As0}) and take in account the relations $s_k=s_{2m-k}$).

The solution of the equations for $G_k$ has the form:
\bee\label{AGK}
G_k(t)= \frac {e^{it}f_k(t)} {(e^{it}-\lambda^k)^2},
\eee
where
\bee\label{A1.47}
\!\!\!\!\!\!\!\!\!\! f_k(t)&=&\frac{2\lambda ^{k}}{m}X_{+}[1-\cos (t\mu _{0})]+(-1)^{k}\frac{2\lambda
^{k}}{m}X_{-}[1-\cos (t\mu _{1})]+(1-\lambda ^{k})^{2}s_{k}+
\nn
&+&\frac{2i}{m}(e^{it}-\lambda ^{k})[\mu _{0}X_{+}\sin (t\mu
_{0})+(-1)^{k}\mu _{1}X_{-}\sin (t\mu _{1})],
\eee
and where $X_{\pm}:=X_1 \pm X_2$.

The following proposition is analogous to Proposition \ref{1.7.8}
but its proof is slightly more difficult:

\proposition\label{A1.7.8}
{\it The trace on the algebra \HG \ is degenerate  if and only if the
generating functions $F_p^{tr}$ defined by  formula (\ref{A1.33})
have the following form
\bee\label{Awid1}
F_p^{tr}(t)=\sum_{j=1}^{j_p} \exp(t \omega_{j,p})\varphi_{j,p}(t),
\eee
where $\omega_{j,p}\in \mathbb C$ and $\varphi_{j,p}\in \mathbb C [t]$.
}


\subsection{The degeneracy conditions for the trace}\label{Adegcond}

We now find the values of the parameters $\mu_0$ and $\mu_1$ for which there exists
a nonzero trace $tr$, (i.e., the values
$s_k$ (\ref{A1.2.3}), not all zero) such that the generating functions (\ref{AGK})
are of the form (\ref{wid1}).
Obviously, it is necessary that the numerator of Eq. (\ref{AGK}) contains all zeros
of the denominator of this expression.
The denominator of $G_k$ vanishes at the points
\be     \nonumber
t_{k,l}=\frac {\pi} m k+ 2\pi l, \text{ where }
\quad l=0,\, \pm 1,\, \pm 2, \, ...
\ee

It so happens that it is sufficient to consider only
the points $t_{k,0}$.

Set
\be     \nonumber
s^{\prime}_k := s_{k} \sin^2\left(\frac {\pi k} {2m} \right),
\quad k=1,...,2m-1, \qquad s^\prime_0=0.
\ee

Then the system of linear equations for $s^\prime_k$
has the form
\bee\label{Aur1}
&& \!\!\!\!\!\!\!\!\!\!\!\!\!\!\!\!\!\!\!\!
\left(1 - \cos \left(\frac {\pi}{ m} k\mu _{0}\right)\right)X_{+}
+ (-1)^{k} \left(1 - \cos \left(\frac {\pi} {m} k\mu _{1}\right)\right)X_{-}
=2 m  s^\prime_{k},  \quad k=1,...,2m-1,\\
&&  \!\!\!\!\!\!\!\!\!\!\!\!\!\!\!\!\!\!\!\!
s^\prime_{2m-r}=s^\prime_r, \quad r=1,...,m, \label{Aur2} \\
\label{Aur3}
&& \!\!\!\!\!\!\!\!\!\!\!\!\!\!\!\!\!\!\!\!
X_{\pm}=X_1 \pm X_2,\\
\label{Aur4}
&& \!\!\!\!\!\!\!\!\!\!\!\!\!\!\!\!\!\!\!\!
 X_1= \sum_{1\le l\le m-1} s^\prime_{2l},               \\
\label{Aur5}
&&  \!\!\!\!\!\!\!\!\!\!\!\!\!\!\!\!\!\!\!\!
X_2= \sum_{0\le l\le m-1}s^\prime_{2l+1},
\eee
and the parameters $\mu_0$ and $\mu_1$ are defined from the condition
that this system has a nonzero solution.

Eqs. (\ref{Aur1}) -- (\ref{Aur5}) imply
that the dimension of the space of solutions is $\leqslant$2
and we can take the values $X_1$ and $X_2$ as parameters determining
the solutions.


\theorem\label{Ath1}
{\it Let $m \geqslant 2$. Then
the system of equations (\ref{Aur1})-(\ref{Aur5}) has nonzero solutions
at the following values of the parameters $\mu_0$ and $\mu_1$ only:
\bee\label{Ak1}
&& \mu_0 \in \mathbb Z\diagdown m\mathbb Z, \qquad \mu_1 \in \mathbb Z\diagdown m\mathbb Z,
\\
\label{Ak3}
&& \mu_0 \in \mathbb Z\diagdown m\mathbb Z, \qquad \text{any } \mu_1,
\\
\label{Ak4}
&& \mu_1 \in \mathbb Z\diagdown m\mathbb Z, \qquad \text{any } \mu_0,
\\
\label{Ak5}
&& \mu_0=\pm \mu_1+ m(2l+1), \qquad l=0,\,\pm 1,\, \pm 2,\,...
\eee
Here,

1. In  case (\ref{Ak1}), the system of equations (\ref{Aur1})-(\ref{Aur5}) has a two-parameter family
of solutions;

2.
In case (\ref{Ak3}), if $\mu_1 \notin \mathbb Z \diagdown m\mathbb Z$,
then the system of equations (\ref{Aur1})-(\ref{Aur5}) has a one-parameter family of solutions
with $X_-=0$,

3.
In case (\ref{Ak4}), if $\mu_0 \notin \mathbb Z \diagdown m\mathbb Z$,
then the system of equations (\ref{Aur1})-(\ref{Aur5}) has a one-parameter family of solutions
with $X_+=0$,

4.
In case (\ref{Ak5}), if $\mu_0, \mu_1 \notin \mathbb Z \diagdown m\mathbb Z$,
then the system of equations (\ref{Aur1})-(\ref{Aur5}) has a one-parameter family of solutions
with $X_1=0$.
}

\remark
Theorem
A.5
is proved for $m\geqslant 2$,
nevertheless it describes also the case $m=1$ correctly.

If $m=1$ then the cases (\ref{Ak1}) -- (\ref{Ak4})  disappear,
and the case (\ref{Ak5}) shows that
\be\label{Arem}
\text{at least one of $\nu_1$ and $\nu_2$
is half-integer.}
\ee
Because
$H_{1,\nu_1, \nu_2}(I_2(2)) \simeq H_{1,\nu_1}(A_1)\otimes H_{1,\nu_2}(A_1)$,
the statement (\ref{Arem}) follows also from \cite{V}, where
the singular values of $\nu$ and ideals in $H_{1,\nu}(A_1)$
were found.


\end{document}